\documentclass[12pt,a4paper]{amsart}
\usepackage[T1]{fontenc}
\usepackage{lmodern}
\usepackage{amssymb}
\usepackage[noadjust]{cite}
\usepackage{booktabs}
\usepackage{xurl}
\usepackage{hyperref}

\theoremstyle{plain}
\newtheorem{theorem}{Theorem}[section]
\newtheorem{corollary}[theorem]{Corollary}
\newtheorem{lemma}[theorem]{Lemma}
\newtheorem{proposition}[theorem]{Proposition}
\newtheorem{conjecture}{Conjecture}[section]

\theoremstyle{definition}

\theoremstyle{remark}
\newtheorem{remark}[theorem]{Remark}
\newtheorem{example}[theorem]{Example}
\numberwithin{equation}{section}

\DeclareMathOperator{\inv}{inv}
\DeclareMathOperator{\per}{per}
\DeclareMathOperator{\diag}{diag}

\hypersetup{colorlinks=true,linkcolor=blue,citecolor=blue,urlcolor=blue,
 pdftitle={Real Positive Definite Counterexamples to a Monotonicity Conjecture for the q-Permanent},
 pdfauthor={Carlos M. da Fonseca and Dev D. Goyal}}

\begin{document}

\title[Counterexamples to a $q$-permanent monotonicity]{Real positive definite counterexamples to a monotonicity conjecture for the $q$-permanent}

\author[C.M. da Fonseca]{Carlos M. da Fonseca}
\address[C.M. da Fonseca]{Kuwait College of Science and Technology, Doha District, Safat 13133, Kuwait}
\email{c.dafonseca@kcst.edu.kw}

\author[D.D. Goyal]{Dev D. Goyal}
\address[D.D. Goyal]{Independent Researcher, New Delhi, India}
\email{dvgyl@outlook.com}

\subjclass[2020]{Primary 15A15; Secondary 15B48, 05A05, 68V20}
\keywords{$q$-permanent, $\mu$-permanent, positive definite matrix, Bapat's conjecture, monotonicity, inversions, formal verification}

\date{November 9, 2026}

\dedicatory{To the memory of Ravindra B. Bapat}

\begin{abstract}
The $q$-permanent $P_q$ of an $n\times n$ matrix weights the term of each permutation $\sigma$ by $q^{\inv(\sigma)}$, where $\inv(\sigma)$ is the number of inversions of $\sigma$. It interpolates between the determinant ($q=-1$) and the permanent ($q=1$). Extending a conjecture of Bapat on $[-1,1]$, it was conjectured in 2010 that the $q$-permanent of every non-diagonal real symmetric positive definite matrix is strictly increasing on an interval $(\varepsilon,+\infty)$ with $\varepsilon<-1$. In this note, we disprove this conjecture. An explicit real symmetric positive definite integer matrix of order five, whose $q$-permanent has only nine nonzero terms, satisfies $P_{14}(A)>P_{15}(A)$. A sparse family shows that for every real $r>1$ some real symmetric positive definite matrix has a $q$-permanent with negative derivative at $q=r$. The counterexamples can be chosen to be correlation matrices without zero entries or, separately, to have integer entries. We also show that the smallest order of a counterexample is four, and that matrices supported on three coordinates separated by identity blocks fail much closer to $q=1$: in order five already near $q=1.7$, and in order $O(1/(r-1))$ at any prescribed $q=r>1$.
\end{abstract}

\maketitle

\nocite{andelic2018,bapat1992,bapatlal1994,bozejko1991,briand2005,fonseca2005,fonseca2010,fonseca2019,hashimoto1992,kamkuo2025,kitamura,lal1992,lal1998,marcus1963,mitchell,noumi1989,noumi1993,tagawa1991,tagawa1993,yang1991}

\section{Introduction}

Let $A=(a_{ij})$ be an $n\times n$ complex matrix and let $q$ be a real number. The \emph{$q$-permanent} of $A$ is defined by
\begin{equation}\label{eq:def}
 P_q(A)=\sum_{\sigma\in S_n}q^{\inv(\sigma)}\prod_{i=1}^n a_{i\sigma(i)},
 \qquad
 \inv(\sigma)=\#\{(i,j):\ i<j,\ \sigma(i)>\sigma(j)\},
\end{equation}
with the convention $0^0=1$. Thus $P_{-1}(A)=\det A$, $P_1(A)=\per A$, and $P_0(A)=a_{11}\cdots a_{nn}$. Unlike the determinant and the permanent, $P_q(A)$ depends on the order of the rows and columns: a simultaneous permutation of rows and columns may change it.  Throughout, matrices are taken in the displayed order, and $P'_q(A)$ denotes the derivative of $P_q(A)$ with respect to~$q$.

The $q$-permanent was introduced independently, under various names, in matrix theory and combinatorics \cite{bapat1992,bapatlal1994,yang1991} and in the study of quantum groups and $q$-deformed commutation relations \cite{bozejko1991,hashimoto1992,noumi1989,noumi1993,tagawa1991,tagawa1993}. In part of the literature it is called the $\mu$-permanent and denoted by $P_\mu(A)$ (see \cite{fonseca2005,fonseca2010,andelic2018,fonseca2019} and the references therein). Related interpolations between the determinant and the permanent are considered in \cite{briand2005}. More recently, Kam and Kuo \cite{kamkuo2025} suggested that $q$-permanents with $q$ close to $1$ could describe off-site crosstalk in boson sampling with transmon qubits.

Bo\.zejko and Speicher \cite{bozejko1991} proved that the function $\sigma\mapsto q^{\inv(\sigma)}$ is positive definite on $S_n$ for $q\in[-1,1]$. Using this, Bapat \cite{bapat1992} showed that $P_q(A)\geqslant0$, for every Hermitian positive semidefinite matrix $A$ and every $q\in[-1,1]$, and proposed the following conjecture, which he proved for $n\leqslant3$ (see also \cite{bapatlal1994,lal1998} for further inequalities for the $q$-permanent).

\begin{conjecture}\label{conj:BL}
If $A$ is a non-diagonal Hermitian positive definite matrix, then $P_q(A)$ is a strictly increasing function of $q\in[-1,1]$.
\end{conjecture}

By continuity, Conjecture~\ref{conj:BL} would imply that $\det A\leqslant P_q(A)\leqslant\per A$ for every Hermitian positive semidefinite matrix $A$ and every $q\in[-1,1]$. Since $P_0(A)=a_{11}\cdots a_{nn}$, this would generalize both Hadamard's inequality $\det A\leqslant a_{11}\cdots a_{nn}$ and its permanental analogue $\per A\geqslant a_{11}\cdots a_{nn}$, due to Marcus \cite{marcus1963}.

Conjecture~\ref{conj:BL} was proved for irreducible tridiagonal positive definite matrices by Lal \cite{lal1992} and, independently, in \cite{fonseca2005}, where it is deduced from a more general result for positive definite matrices whose graph is a tree. That result requires, in addition, a noncrossing condition on the labelling of the vertices of the tree, which is missing from its statement in \cite{fonseca2005} (see also \cite{andelic2018}). Mitchell \cite{mitchell} proved that Conjecture~\ref{conj:BL} is equivalent to its extension to non-diagonal singular positive semidefinite matrices without zero rows and established this extension for matrices of rank one and for matrices of order three. A Lean-verified complex Hermitian counterexample to Conjecture~\ref{conj:BL} of order $144$ has recently been made public by Kitamura \cite{kitamura}.

In \cite{fonseca2010}, where $A>0$ denotes a real symmetric positive definite matrix, the first author proposed to extend Conjecture~\ref{conj:BL} beyond $[-1,1]$. In the form restated in \cite[Conjecture~2]{fonseca2019}, the extension reads as follows.

\begin{conjecture}[{\cite{fonseca2010}, \cite[Conjecture~2]{fonseca2019}}]\label{conj:ext}
For every non-diagonal real symmetric positive definite matrix $A$ there exists $\varepsilon<-1$ such that $P_q(A)$ is a strictly increasing function of $q\in(\varepsilon,+\infty)$.
\end{conjecture}

The word \emph{non-diagonal} is implicit in \cite{fonseca2019}, since $P_q(A)$ does not depend on $q$ when $A$ is diagonal. The formulation in \cite[Section~2, Conjecture~1]{fonseca2010} asks for strict monotonicity on an interval whose left endpoint is at most $-1$ and whose right endpoint is $+\infty$. Either formulation implies that
\begin{equation}\label{eq:implied}
 P_{q_1}(A)<P_{q_2}(A),\qquad\text{whenever } -1<q_1<q_2,
\end{equation}
for every non-diagonal real symmetric positive definite matrix $A$.

Since $(\varepsilon,+\infty)\supset[-1,1]$, Conjecture~\ref{conj:ext} strengthens the real symmetric case of Conjecture~\ref{conj:BL}. Kitamura's matrix has nonreal entries, so it is not a counterexample to Conjecture~\ref{conj:ext}, and we know of no general procedure that turns a complex Hermitian counterexample into a real symmetric one. Kitamura's complex example does not by itself decide the real symmetric case of Conjecture~\ref{conj:BL}. The constructions in this note concern $q>1$ and do not themselves settle monotonicity on $[-1,1]$.

In this note, we show that \eqref{eq:implied} fails, so that Conjecture~\ref{conj:ext} is false in each of its forms, already for real symmetric matrices with integer entries. This will be done in Sections~\ref{sec:five}--\ref{sec:allq}. In Section~\ref{sec:small} we show that no counterexample of order at most three exists and give one of order four, so that the smallest order of a counterexample is four. Next, we study \emph{gapped} matrices, supported on three coordinates separated by identity blocks. Their $q$-permanent has at most six terms. They give an integer positive definite matrix of order five with $P_2(A)>P_3(A)$ and even $P_3(A)<0$, failure at every $q=r>1$ in order $O(1/(r-1))$, and explicit upper bounds for the smallest failing parameter in each order $4\leqslant n\leqslant9$.

The main counterexample theorems, the small-order minimality result, the gapped-family identity and quantitative bounds, and all six witnesses in Table~\ref{tab:thresholds} have been formalized in Lean~4 (version~4.35.0, with Mathlib at commit \texttt{56a0058}). The Sturm root isolations, displayed principal minors, and finite lists of sufficient dimensions are checked independently in exact integer and rational arithmetic. The supplementary material contains the formal proofs, exact certificates and reproducible computations. In all of our explicit examples, $P_q(A)$ is strictly increasing on $[-1,1]$.

\section{A counterexample of order five}\label{sec:five}

We begin with a single explicit counterexample to Conjecture~\ref{conj:ext}. It is a sparse integer matrix of order five. Only nine permutations contribute to its $q$-permanent, so $P_q(A)$ can be computed by hand, and the failure of \eqref{eq:implied} comes down to comparing two integers.

\begin{theorem}\label{thm:five}
The matrix
\[
 A=\begin{pmatrix}
 5&0&0&4&1\\
 0&5&0&-1&4\\
 0&0&1&0&0\\
 4&-1&0&5&0\\
 1&4&0&0&5
 \end{pmatrix}
\]
is positive definite, and
\begin{equation}\label{eq:five-poly}
 P_q(A)=625+25q^3+800q^5+25q^7+256q^8-32q^9+q^{10}.
\end{equation}
In particular,
\[
 P_{14}(A)=8\,968\,822\,249>7\,441\,960\,000=P_{15}(A),
\]
so $A$ violates \eqref{eq:implied} and Conjecture~\ref{conj:ext} is false.
\end{theorem}

\begin{proof}
Reordering the coordinates as $(1,2,4,5,3)$ transforms $A$ into
\[
 \begin{pmatrix}5I_2&R\\R^{T}&5I_2\end{pmatrix}\oplus[1],\qquad
 R=\begin{pmatrix}4&1\\-1&4\end{pmatrix},\qquad RR^{T}=17I_2 .
\]
Hence the eigenvalues of $A$ are $1$ and $5\pm\sqrt{17}$, each of the latter with multiplicity two, and $A$ is positive definite. Equivalently, the leading principal minors of $A$ are $5$, $25$, $25$, $40$ and $64$.

For the $q$-permanent we return to the original order. Row and column~$3$ of $A$ vanish off the diagonal, so every permutation $\sigma$ with a nonzero term fixes~$3$. The remaining nonzero off-diagonal entries lie in the positions $\{1,2\}\times\{4,5\}$ and $\{4,5\}\times\{1,2\}$. Consequently, the restriction of $\sigma$ to $\{1,2,4,5\}$ is the identity, one of the four transpositions $(1\ 4)$, $(1\ 5)$, $(2\ 4)$, $(2\ 5)$, one of the two products $(1\ 4)(2\ 5)$, $(1\ 5)(2\ 4)$, or one of the two $4$-cycles $(1\ 4\ 2\ 5)$, $(1\ 5\ 2\ 4)$. These nine terms are listed in Table~\ref{tab:terms}. Summing them gives \eqref{eq:five-poly}, and the values at $q=14$ and $q=15$ follow by substitution.
\end{proof}

\begin{table}[htbp]
\centering
\begin{tabular}{@{}lcr@{}}
\toprule
Permutation (one-line notation) & $\inv(\sigma)$ & $\prod_i a_{i\sigma(i)}$\\
\midrule
$(1,2,3,4,5)$ & $0$ & $625$\\
$(1,4,3,2,5)$ & $3$ & $25$\\
$(1,5,3,4,2)$ & $5$ & $400$\\
$(4,2,3,1,5)$ & $5$ & $400$\\
$(5,2,3,4,1)$ & $7$ & $25$\\
$(4,5,3,1,2)$ & $8$ & $256$\\
$(4,5,3,2,1)$ & $9$ & $-16$\\
$(5,4,3,1,2)$ & $9$ & $-16$\\
$(5,4,3,2,1)$ & $10$ & $1$\\
\bottomrule
\end{tabular}
\caption{The nine nonzero terms of $P_q(A)$ for the matrix $A$ of Theorem~\ref{thm:five}.}
\label{tab:terms}
\end{table}

We observe that, by \eqref{eq:five-poly}, $P'_q(A)=q^2g(q)$ with
\[
 g(q)=75+4000q^2+175q^4+2048q^5-288q^6+10q^7 .
\]
A Sturm sequence computation shows that $g$ is positive on $[-1,\alpha)$, negative on $(\alpha,\beta)$ and positive on $(\beta,+\infty)$, where $13.35<\alpha<13.36$ and $15.55<\beta<15.56$. Hence $P_q(A)$ is strictly increasing on $[-1,\alpha]$ and on $[\beta,+\infty)$, and strictly decreasing on $[\alpha,\beta]$. For instance, $P'_{15}(A)=-1\,563\,030\,000$. In particular, $P_q(A)$ is strictly increasing on $[-1,1]$, in agreement with Conjecture~\ref{conj:BL}.

\section{A sparse family}\label{sec:family}

Let $m\geqslant0$ be an integer and $n=m+4$. We call
\[
 i=1,\qquad j=2,\qquad k=m+3,\qquad \ell=m+4
\]
the \emph{active} coordinates, and the remaining $m$ coordinates $3,\dots,m+2$ \emph{passive}. For real numbers $a\geqslant0$ and $d>\sqrt{a^2+1}$, let $A_m(a,d)$ be the real symmetric matrix of order $n$ that coincides with the identity on the passive coordinates, has zero entries between passive and active coordinates, and whose principal submatrix on the active coordinates, in increasing order, is
\begin{equation}\label{eq:core}
 \begin{pmatrix}
 d&0&a&1\\
 0&d&-1&a\\
 a&-1&d&0\\
 1&a&0&d
 \end{pmatrix}.
\end{equation}
As in the proof of Theorem~\ref{thm:five}, this block equals $\left(\begin{smallmatrix}dI_2&R\\R^{T}&dI_2\end{smallmatrix}\right)$ with $R=\left(\begin{smallmatrix}a&1\\-1&a\end{smallmatrix}\right)$ and $RR^{T}=(a^2+1)I_2$. Its eigenvalues are therefore $d\pm\sqrt{a^2+1}$, each with multiplicity two, and $A_m(a,d)$ is positive definite. We stress that the passive coordinates lie between $j$ and $k$ and they cannot be deleted when inversions are counted.

\begin{lemma}\label{lem:family}
For every non-negative integer $m$,
\begin{equation}\label{eq:family}
 P_q(A_m(a,d))
 =d^4+d^2q^{2m+1}\bigl(1+2a^2q^2+q^4\bigr)+q^{4m+4}\bigl(a^2-q\bigr)^2.
\end{equation}
\end{lemma}

\begin{proof}
Every permutation $\sigma$ with a nonzero term fixes the passive coordinates. As in the proof of Theorem~\ref{thm:five}, its restriction to $\{i,j,k,\ell\}$ is one of the nine permutations in Table~\ref{tab:family}, written in cycle notation.

For a transposition of the positions $r<s$, the number of inversions is $2(s-r)-1$. This gives the entries $2m+1$, $2m+3$, and $2m+5$ in the table. Each of the remaining four permutations maps $\{i,j\}$ onto $\{k,\ell\}$ and $\{k,\ell\}$ onto $\{i,j\}$. Hence every passive coordinate forms an inversion with each of the four active coordinates, contributing $4m$ inversions in total, while the pairs of active coordinates contribute $4$, $5$, $5$ and $6$ inversions, respectively. Summing the entries of Table~\ref{tab:family} gives
\[
 d^4+d^2q^{2m+1}+2a^2d^2q^{2m+3}+d^2q^{2m+5}
 +a^4q^{4m+4}-2a^2q^{4m+5}+q^{4m+6},
\]
which is \eqref{eq:family}.
\end{proof}

\begin{table}[htbp]
\centering
\begin{tabular}{@{}lcc@{}}
\toprule
Restriction of $\sigma$ to $\{i,j,k,\ell\}$ & $\inv(\sigma)$ & $\prod_i a_{i\sigma(i)}$\\
\midrule
identity & $0$ & $d^4$\\
$(j\ k)$ & $2m+1$ & $d^2$\\
$(i\ k)$,\ $(j\ \ell)$ & $2m+3$ & $a^2d^2$ each\\
$(i\ \ell)$ & $2m+5$ & $d^2$\\
$(i\ k)(j\ \ell)$ & $4m+4$ & $a^4$\\
$(i\ k\ j\ \ell)$,\ $(i\ \ell\ j\ k)$ & $4m+5$ & $-a^2$ each\\
$(i\ \ell)(j\ k)$ & $4m+6$ & $1$\\
\bottomrule
\end{tabular}
\caption{The potentially nonzero terms of $P_q(A_m(a,d))$.}
\label{tab:family}
\end{table}

Note that the parameter $m$ does not affect the nonzero entries of the active block: it only shifts the exponents of $q$ in \eqref{eq:family}. This is what allows the negative contribution of the last summand to be made dominant in Section~\ref{sec:allq}. The matrix of Theorem~\ref{thm:five} is $A_1(4,5)$, with passive coordinate~$3$, and \eqref{eq:five-poly} is the case $m=1$, $a=4$, $d=5$ of \eqref{eq:family}. Table~\ref{tab:terms} is the corresponding specialization of Table~\ref{tab:family}.

\section{Failure at every prescribed parameter above one}\label{sec:allq}

We start this section by observing that the last summand of \eqref{eq:family} can have a negative derivative. Choosing $a$ and $d$ as functions of $m$ and of a fixed $r>1$, we make this negative contribution at $q=r$ grow faster with $m$ than the positive derivative of the middle summand.

\begin{theorem}\label{thm:allq}
For every real $r>1$ there exist an integer $n\geqslant5$ and a real symmetric positive definite matrix $A$ of order $n$ such that $P'_r(A)<0$. Consequently, for every $\delta>0$ there is a real symmetric positive definite matrix $A$ such that $P_q(A)$ is not increasing on $[1,1+\delta]$.
\end{theorem}

\begin{proof}
Fix $r>1$ and choose an integer $m\geqslant1$ such that
\begin{equation}\label{eq:m-bound}
 r^{2m+5}>(4m+4)(2m+5)\Bigl(\frac{5r}{4}+2\Bigr)\Bigl(1+\frac52r^3+r^4\Bigr).
\end{equation}
Such an $m$ exists, since the left-hand side grows exponentially in $m$, while the right-hand side is a quadratic polynomial in $m$ whose coefficients depend only on $r$. Set
\[
 K=4m+4,\qquad H=2m+1,\qquad
 t=r\Bigl(1+\frac1K\Bigr),\qquad
 a=\sqrt{t},\qquad d=\sqrt{t+2}.
\]
Then $d^2=a^2+2>a^2+1$, so $A=A_m(a,d)$ is positive definite of order $n=m+4\geqslant5$. By Lemma~\ref{lem:family},
\[
 P_q(A)=(t+2)^2+(t+2)\,q^H\bigl(1+2tq^2+q^4\bigr)+q^K(t-q)^2,
\]
where $t$ does not depend on $q$. Since $t-r=r/K$, the derivative of the last summand at $q=r$ equals
\[
 Kr^{K-1}(t-r)^2-2r^K(t-r)
 =\frac{r^{K+1}}{K}-\frac{2r^{K+1}}{K}
 =-\frac{r^{K+1}}{K}.
\]
The derivative of the middle summand at $q=r$ equals
\[
 (t+2)\,r^{H-1}\bigl(H+2t(H+2)r^2+(H+4)r^4\bigr).
\]
Since $K\geqslant8$, we have $t\leqslant9r/8<5r/4$, so that $t+2\leqslant\frac{5r}4+2$ and $2tr^2\leqslant\frac52r^3$. Moreover, $H<H+2<H+4=2m+5$. Hence
\begin{align*}
 (t+2)\,r^{H-1}\bigl(H+2t(H+2)r^2+(H+4)r^4\bigr)
 &\leqslant\Bigl(\frac{5r}{4}+2\Bigr)(2m+5)\,r^{H-1}\Bigl(1+\frac52r^3+r^4\Bigr)\\
 &<\frac{r^{H-1}\,r^{2m+5}}{K}=\frac{r^{K+1}}{K},
\end{align*}
where the strict inequality is \eqref{eq:m-bound} multiplied by $r^{H-1}/K$, and the last equality uses $(H-1)+(2m+5)=K+1$. Therefore $P'_r(A)<0$.

For the last assertion, given $\delta>0$, apply the first part with $r=1+\delta/2$. Since $P'_r(A)<0$, the polynomial $P_q(A)$ is strictly decreasing on a neighbourhood of $r$ contained in $[1,1+\delta]$.
\end{proof}

\begin{remark}\label{rem:orders}
The order of the matrices in the proof grows as $r$ decreases to $1$. For $r=101/100$, $11/10$, $3/2$, $2$, and $15$, the least $m\geqslant1$ satisfying \eqref{eq:m-bound} is $925$, $69$, $12$, $6$, and $1$, respectively.
\end{remark}

\begin{corollary}\label{cor:correlations}
The matrix of Theorem~\ref{thm:five} and the matrices of Theorem~\ref{thm:allq} may be replaced by real positive definite correlation matrices with the same properties. Moreover, for every $r>1$ there is a real positive definite correlation matrix $A$ without zero entries such that $P'_r(A)<0$.
\end{corollary}

\begin{proof}
If $D=\diag(s_1,\dots,s_n)$ with $s_1,\dots,s_n>0$, then $(DAD)_{i\sigma(i)}=s_is_{\sigma(i)}a_{i\sigma(i)}$, and therefore, term by term,
\begin{equation}\label{eq:scaling}
 P_q(DAD)=(\det D)^2\,P_q(A).
\end{equation}
Choosing $s_i=a_{ii}^{-1/2}$ normalizes the diagonal. Positive definiteness is preserved and the factor $(\det D)^2>0$ does not depend on $q$. For the matrix of Theorem~\ref{thm:five}, this gives a correlation matrix $C$ whose off-diagonal entries lie in $\{0,\pm\frac15,\frac45\}$ and $P_q(C)=P_q(A)/625$.

For the second assertion, let $A$ be a matrix given by Theorem~\ref{thm:allq} and let $\mathbf{1}$ be the all-ones vector. For every $\varepsilon>0$, the matrix $A+\varepsilon\mathbf{1}\mathbf{1}^{T}$ is positive definite, and for sufficiently small $\varepsilon>0$ all its entries are nonzero. Since $P'_r(\cdot)$ is a polynomial in the entries, it remains negative for sufficiently small $\varepsilon$. Normalizing the diagonal as above completes the proof.
\end{proof}

\begin{corollary}\label{cor:integer}
For every real $r>1$ there is a real symmetric positive definite matrix $A$ with integer entries such that $P'_r(A)<0$.
\end{corollary}

\begin{proof}
Positive definiteness and the inequality $P'_r(A)<0$ are open conditions on the real vector space of symmetric matrices of order $n$. Hence the matrix of Theorem~\ref{thm:allq} can be approximated by a rational symmetric matrix $B$ satisfying both conditions. If $c$ is a positive common denominator of the entries of $B$, then $cB$ has integer entries, is positive definite, and $P'_r(cB)=c^nP'_r(B)<0$.
\end{proof}

\section{Small orders}\label{sec:small}

The family of Section~\ref{sec:family} has order at least four; the construction used in Theorem~\ref{thm:allq} has order at least five. In this section, with exact integer and rational arithmetic, we show that order four already suffices and that order three does not.

For $q\in[-1,1]$, the case $n\leqslant3$ of Conjecture~\ref{conj:BL} is due to Bapat \cite{bapat1992}. In the real symmetric case, the following elementary proposition extends it to all real~$q$.

\begin{proposition}\label{prop:small}
If $A$ is a non-diagonal real symmetric positive definite matrix of order $n\leqslant3$, then $P_q(A)$ is strictly increasing on $\mathbb{R}$.
\end{proposition}

\begin{proof}
By \eqref{eq:scaling} we may assume that $A$ is a correlation matrix, so that its off-diagonal entries have absolute value less than one. For $n=2$, $P_q(A)=1+qa_{12}^2$, with $a_{12}\ne0$. For $n=3$, write $x=a_{12}$, $y=a_{23}$ and $z=a_{13}$. Then
\[
 P_q(A)=1+(x^2+y^2)q+2xyz\,q^2+z^2q^3,\qquad
 P'_q(A)=x^2+y^2+4xyz\,q+3z^2q^2 .
\]
If $z=0$, the derivative equals $x^2+y^2>0$. If $z\ne0$, minimizing the quadratic in $q$ gives
\[
 P'_q(A)\geqslant x^2+y^2-\tfrac43x^2y^2
 =x^2\bigl(1-\tfrac23y^2\bigr)+y^2\bigl(1-\tfrac23x^2\bigr),
\]
which is positive unless $x=y=0$. In such a case the derivative is $3z^2q^2$, which vanishes only at $q=0$. In all cases $P_q(A)$ is strictly increasing on $\mathbb{R}$.
\end{proof}

\begin{theorem}\label{thm:four}
The matrix
\[
 B=\begin{pmatrix}
 52&29&48&-1\\
 29&52&38&42\\
 48&38&52&14\\
 -1&42&14&52
 \end{pmatrix}
\]
is positive definite and
\begin{multline*}
 P_q(B)=7\,311\,616+6\,708\,624\,q+7\,989\,796\,q^2+12\,606\,008\,q^3\\
 +3\,867\,696\,q^4-150\,512\,q^5+1\,444\,q^6 .
\end{multline*}
In particular,
\[
 P_{44}(B)=1\,242\,343\,329\,792>1\,242\,122\,338\,096=P_{45}(B).
\]
\end{theorem}

\begin{proof}
The leading principal minors of $B$ are $52$, $1863$, $7772$ and $6432$, so $B$ is positive definite. The polynomial is obtained by expanding \eqref{eq:def} over the $24$ permutations of $S_4$, and the two values follow by substitution.
\end{proof}

A Sturm sequence computation shows that $P_q(B)$ is strictly increasing on $[-1,\alpha']$ and on $[\beta',+\infty)$ and strictly decreasing on $[\alpha',\beta']$, where $43.95<\alpha'<43.96$ and $45.15<\beta'<45.16$. The following observation explains the negative coefficient in this example. For a matrix of order four, the coefficients of $q^6$ and $q^5$ in $P_q$ are
\[
 c_6=a_{14}^2a_{23}^2,\qquad c_5=a_{14}^2a_{22}a_{33}+2a_{14}a_{13}a_{23}a_{24},
\]
and $P'_q=c_1+2c_2q+3c_3q^2+q^3\bigl(4c_4+5c_5q+6c_6q^2\bigr)$, where $c_h$ denotes the coefficient of $q^h$. Since $c_5=a_{14}\bigl(a_{14}a_{22}a_{33}+2a_{13}a_{23}a_{24}\bigr)$, we have $c_5<0$ whenever
\[
 a_{13}a_{14}a_{23}a_{24}<0\qquad\text{and}\qquad 0<|a_{14}|<\frac{2|a_{13}a_{23}a_{24}|}{a_{22}a_{33}}.
\]
If, in addition, $a_{23}\neq0$, then $c_6>0$, and the quadratic factor $4c_4+5c_5q+6c_6q^2$ can become negative for $q$ of the order of $|c_5|/c_6$. For $B$ this happens near $q\approx44$.

\begin{corollary}\label{cor:minimal}
The smallest order of a non-diagonal real symmetric positive definite matrix violating \eqref{eq:implied} is four.
\end{corollary}

\begin{proof}
Combine Proposition~\ref{prop:small} and Theorem~\ref{thm:four}.
\end{proof}

\section{Gapped matrices and failure near \texorpdfstring{$q=1$}{q=1}}\label{sec:gapped}

The counterexamples of Sections~\ref{sec:five} and~\ref{sec:small} fail only for fairly large values of~$q$. A different and even simpler construction shows that matrices of small order already fail close to $q=1$.

Let $n\geqslant3$ and $1<p<n$. For a real symmetric matrix $C=(c_{st})$ of order three, indexed by $\{1,p,n\}$, let $G_{n,p}(C)$ be the real symmetric matrix of order $n$ whose principal submatrix on the coordinates $1,p,n$ is $C$, which coincides with the identity on the remaining coordinates, and whose other entries are zero. We call $G_{n,p}(C)$ a \emph{gapped} matrix. Its $m_1=p-2$ coordinates between $1$ and $p$ and its $m_2=n-p-1$ coordinates between $p$ and $n$ are passive, and $G_{n,p}(C)$ is positive definite if and only if $C$ is.

\begin{lemma}\label{lem:gapped}
With $M=m_1+m_2=n-3$,
\begin{multline*}
 P_q\bigl(G_{n,p}(C)\bigr)=c_{11}c_{pp}c_{nn}+c_{1p}^2c_{nn}\,q^{2m_1+1}+c_{pn}^2c_{11}\,q^{2m_2+1}\\
 +2c_{1p}c_{pn}c_{1n}\,q^{2M+2}+c_{1n}^2c_{pp}\,q^{2M+3}.
\end{multline*}
\end{lemma}

\begin{proof}
Every permutation with a nonzero term fixes the passive coordinates, so its restriction to $\{1,p,n\}$ is one of the six permutations of this set. By the rule $2(s-r)-1$ used in the proof of Lemma~\ref{lem:family}, the transpositions $(1\ p)$, $(p\ n)$ and $(1\ n)$ have $2m_1+1$, $2m_2+1$ and $2M+3$ inversions, respectively. For each of the two $3$-cycles, the pairs of active coordinates contribute two inversions, and a direct check shows that every passive coordinate forms exactly two inversions with active coordinates; this gives $2M+2$ inversions. Both $3$-cycles have entry product $c_{1p}c_{pn}c_{1n}$.
\end{proof}

For $n=3$, Lemma~\ref{lem:gapped} reduces to the formula in the proof of Proposition~\ref{prop:small}. With passive coordinates, however, the term $2c_{1p}c_{pn}c_{1n}q^{2M+2}$, which is negative when $c_{1p}c_{pn}c_{1n}<0$, can dominate the derivative.

\begin{example}\label{ex:gapped}
Let $G=G_{5,3}(C)$ with $C=\left(\begin{smallmatrix}11&7&2\\7&11&-7\\2&-7&11\end{smallmatrix}\right)$, that is,
\[
 G=\begin{pmatrix}
 11&0&7&0&2\\
 0&1&0&0&0\\
 7&0&11&0&-7\\
 0&0&0&1&0\\
 2&0&-7&0&11
 \end{pmatrix}.
\]
Its leading principal minors are $11$, $11$, $72$, $72$, and $13$, so $G$ is positive definite, and by Lemma~\ref{lem:gapped},
\[
 P_q(G)=1331+1078\,q^3-196\,q^6+44\,q^7,\qquad
 P'_q(G)=14q^2\bigl(231-84q^3+22q^4\bigr).
\]
A Sturm sequence computation shows that $P_q(G)$ is strictly increasing on $[-1,\alpha]$ and on $[\beta,+\infty)$ and strictly decreasing on $[\alpha,\beta]$, where $1.706<\alpha<1.707$ and $3.591<\beta<3.592$. In particular,
\[
 P_2(G)=3043>-16\,219=P_3(G).
\]
Note that $P_3(G)<0<13=\det G=P_{-1}(G)$. Thus, for $q>1$, the $q$-permanent of a positive definite matrix can even be negative, in contrast with Bapat's result for $q\in[-1,1]$ \cite{bapat1992}. Only six permutations contribute to $P_q(G)$, against nine in Theorem~\ref{thm:five}.
\end{example}

The following proposition gives a second proof of Theorem~\ref{thm:allq}, with explicit matrices of much smaller order.

\begin{proposition}\label{prop:linear}
Let $r>1$ be real and let $m\geqslant1$ be an integer such that
\begin{equation}\label{eq:linear}
 r^{2m}>\frac{2(4m+3)}{2m+1}.
\end{equation}
Set $$t=\frac{1}{r^{2m+1}} \quad \mbox{and}\quad C=\left(\begin{smallmatrix}1&1/2&t\\1/2&1&-1/2\\t&-1/2&1\end{smallmatrix}\right).$$ Then $A=G_{2m+3,m+2}(C)$ is a positive definite correlation matrix of order $2m+3$ and
\[
 P'_r(A)=(4m+3)-\tfrac12(2m+1)r^{2m}<0.
\]
In particular, for every $r>1$, there is a real positive definite correlation matrix $A$ of order at most $2\lceil \ln 5/(2\ln r)\rceil+3$ with $P'_r(A)<0$.
\end{proposition}

\begin{proof}
The leading principal minors of $C$ are $1$, $\frac34$ and $\det C=\frac12-\frac t2-t^2=(\frac12-t)(1+t)$. By \eqref{eq:linear}, $r^{2m}>4$, so $t<r^{-2m}<\frac14$ and $C$ is positive definite. By Lemma~\ref{lem:gapped} with $m_1=m_2=m$,
\[
 P_q(A)=1+\tfrac12q^{2m+1}-\tfrac12t\,q^{4m+2}+t^2q^{4m+3}.
\]
Differentiating and using $tr^{2m+1}=1$, we obtain
\[
 P'_r(A)=\tfrac12(2m+1)r^{2m}-(2m+1)r^{2m}+(4m+3)=(4m+3)-\tfrac12(2m+1)r^{2m},
\]
which is negative by \eqref{eq:linear}. Finally, $2(4m+3)/(2m+1)\leqslant\frac{14}3<5$, for $m\geqslant1$, so \eqref{eq:linear} holds whenever $r^{2m}\geqslant5$, in particular for $m=\lceil\ln5/(2\ln r)\rceil$ (note that $m\geqslant1$).
\end{proof}

Note that, for $r=101/100$, $11/10$, $3/2$, $2$ and $15$, the least $m$ satisfying \eqref{eq:linear} is $70$, $8$, $2$, $2$ and $1$, which gives orders $143$, $19$, $7$, $7$ and $5$, against $929$, $73$, $16$, $10$ and $5$ in Remark~\ref{rem:orders}. Since $\ln r\sim r-1$ as $r\to1$, Proposition~\ref{prop:linear} shows that order $O\bigl(1/(r-1)\bigr)$ suffices for failure at $q=r$.

For $n\geqslant1$, let $\mathcal{P}_n$ denote the set of real symmetric positive definite matrices of order $n$, and let
\[
 r_n=\inf\{r>1:\ P'_r(A)<0 \text{ for some } A\in\mathcal{P}_n\},
\]
with $\inf\emptyset=+\infty$. By Proposition~\ref{prop:small}, $r_n=+\infty$, for $n\leqslant3$. Since appending identity coordinates at the end does not change $P_q$, the sequence $(r_n)$ is non-increasing, and $r_n\to1$ by Theorem~\ref{thm:allq}. More precisely, for $n\geqslant5$ let $m=\lfloor(n-3)/2\rfloor$, so that $m\geqslant1$, $2m+3\leqslant n$ and $2m\geqslant n-4$. Every $r$ with $r^{2m}\geqslant5$ satisfies \eqref{eq:linear}, so Proposition~\ref{prop:linear} and the monotonicity of $(r_n)$ give
\begin{equation}\label{eq:rn-bound}
 r_n\leqslant r_{2m+3}\leqslant5^{1/(2m)}\leqslant5^{1/(n-4)},\qquad n\geqslant5.
\end{equation}
In particular, $r_n-1=O(1/n)$.

Table~\ref{tab:thresholds} lists rigorous upper bounds for $r_n$, for $4\leqslant n\leqslant9$. Each bound is certified by an explicit rational correlation matrix $A$ with $P'_q(A)<0$ at the listed value of~$q$, checked in exact arithmetic. For $n\geqslant5$, the witness is the gapped matrix $A=G_{n,m_1+2}(C)$, where $C$ has unit diagonal and $c_{1p}=x$, $c_{pn}=y$ and $c_{1n}=z$. It is positive definite because $|x|<1$ and $\det C=1-x^2-y^2-z^2+2xyz>0$, and $P'_q(A)$ is given explicitly by Lemma~\ref{lem:gapped}. The table lists $\det C$ and a decimal approximation of $P'_q(A)$. For $n=4$, the witness is $A=W/1000$, where
\[
 W=\begin{pmatrix}
 1000&500&832&18\\
 500&1000&895&-854\\
 832&895&1000&-534\\
 18&-854&-534&1000
 \end{pmatrix}.
\]
The leading principal minors of $W$ are $1000$, $750\,000$, $1\,391\,000$ and $5\,643\,044$, and
\begin{multline*}
 P_q(W)=10^{12}+1\,336\,181\,000\,000\,q+1\,632\,233\,440\,000\,q^2\\
 +1\,792\,359\,212\,000\,q^3+473\,483\,670\,784\,q^4-22\,569\,212\,160\,q^5+259\,532\,100\,q^6 .
\end{multline*}
A Sturm sequence computation shows that $P_q(W)$ is strictly increasing on $[-1,\alpha'']$ and on $[\beta'',+\infty)$ and strictly decreasing on $[\alpha'',\beta'']$, where $36.33<\alpha''<36.34$ and $38.66<\beta''<38.67$. For instance, $P_{36}(W)>P_{37}(W)>P_{38}(W)$.

\begin{table}[htbp]
\centering\small
\begin{tabular}{@{}cccccc@{}}
\toprule
$n$ & $(m_1,m_2)$ & $(x,\,y,\,z)$ & $\det C$ & $q$ & $P'_q(A)$\\
\midrule
$4$ & -- & $A=W/1000$ & -- & $36.34$ & $-0.1847$\\[2pt]
$5$ & $(1,1)$ & $\bigl(\frac{16}{25},\,\frac{16}{25},\,-\frac{9}{50}\bigr)$ & $\frac{59}{62\,500}$ & $1.71$ & $-0.0791$\\[3pt]
$6$ & $(1,2)$ & $\bigl(\frac{37}{50},\,\frac{51}{100},\,-\frac15\bigr)$ & $\frac{67}{50\,000}$ & $1.45$ & $-0.0377$\\[3pt]
$7$ & $(2,2)$ & $\bigl(\frac{31}{50},\,\frac{31}{50},\,-\frac{23}{100}\bigr)$ & $\frac{369}{250\,000}$ & $1.31$ & $-0.1087$\\[3pt]
$8$ & $(2,3)$ & $\bigl(\frac{17}{25},\,\frac{27}{50},\,-\frac{6}{25}\bigr)$ & $\frac{759}{62\,500}$ & $1.25$ & $-0.2951$\\[3pt]
$9$ & $(3,3)$ & $\bigl(\frac{61}{100},\,\frac{61}{100},\,-\frac14\bigr)$ & $\frac{29}{4\,000}$ & $1.20$ & $-0.2766$\\
\bottomrule
\end{tabular}
\caption{Rational correlation matrices $A$ certifying $r_n\leqslant q$. For $n\geqslant5$, $A=G_{n,m_1+2}(C)$ with $C=\left(\begin{smallmatrix}1&x&z\\x&1&y\\z&y&1\end{smallmatrix}\right)$; for $n=4$, $A=W/1000$, with $W$ displayed above. The derivative values are rounded for display; the inequalities are checked in exact rational arithmetic.}
\label{tab:thresholds}
\end{table}

The certified bounds in Table~\ref{tab:thresholds} decrease sharply between orders four and five. Theorem~\ref{thm:allq} and Proposition~\ref{prop:linear} rule out any extension of the monotonicity of $P_q(A)$ beyond $q=1$ that is uniform over all orders, and show that $r_n-1=O(1/n)$. A matching lower bound is not established here.

The correlation and integer-entry versions of our results are separate statements: we do not claim the existence of an integer correlation matrix with the stated properties. Finally, all of the counterexamples in this note concern $q>1$; they do not themselves settle the real symmetric case of Conjecture~\ref{conj:BL} on $[-1,1]$.

\section*{Declaration on the use of AI tools}

GPT-6 (ProReasoning) in ChatGPT was used to find the initial counterexamples. The authors take full responsibility for the content of this paper.

\bibliographystyle{unsrt}
\bibliography{references}

\end{document}